\documentclass[dvips,twoside]{article}
\usepackage[latin1]{inputenc}
\usepackage[T1]{fontenc}
\usepackage{parskip}
\usepackage{amsmath,amsfonts,amssymb,amsxtra}
\usepackage{latexsym}
\usepackage{theorem}
\usepackage{graphicx,psfrag}
\usepackage{verbatim}

\newtheorem{theo}{Theorem}
\newtheorem{prop}{Proposition}
\newtheorem{lemma}{Lemma}

{\theoremstyle{plain}                       
\theorembodyfont{\rmfamily}
}

{\theoremstyle{plain}                       
\theorembodyfont{\rmfamily}
\newtheorem{Example}{Example}}

{\theoremstyle{plain}                       
\theorembodyfont{\rmfamily}
}

\begin{document}

\title{Kolmogorov-Sinai entropy of a generalized Markov shift}

\author{IVAN WERNER\\
   {\small Email: ivan\_werner@pochta.ru}}
\maketitle

\begin{abstract}
 In this paper we calculate Kolmogorov-Sinai entropy $h_M(S)$ of the generalized
 Markov shift associated with a contractive Markov system (CMS) \cite{Wer1} using the coding map constructed in
 \cite{Wer3}. We show that \[h_M(S)=-\sum\limits_{e\in E}\int\limits_{K_{i(e)}} p_e\log
  p_ed\mu\] where $\mu$ is a unique invariant Borel probability
  measure of the CMS.
\end{abstract}

\section{Introduction}
In  \cite{Wer1} we introduced a theory of {\it contractive Markov
systems (CMS)} which provides a unifying framework in so-called
"fractal" geometry. It extends the known theory of {\it iterated
function systems (IFS)} with place dependent probabilities
\cite{BDEG}\cite{Elton} in a way that it also covers {\it graph
directed constructions} of "fractal" sets \cite{MW}. In
particular, Markov chains associated with such systems also cover
finite Markov chains.

By a {\it Markov system}  we mean a family
\[\left(K_{i(e)},w_e,p_e\right)_{e\in E}\]
(see Fig. 1) where $E$ is the set of edges of a finite directed
(multi)graph $(V,E,i,t)$ ($V:=\{1,...,N\}$ is the set of vertices
of the directed (multi)graph (we do not exclude the case $N=1$),
$i:E\longrightarrow V$ is a map indicating the initial vertex of
each edge and $t:E\longrightarrow V$ is a map indicating the
terminal vertex of each edge), $K_1,K_2,...,K_N$ is a partition of
a metric space $(K,d)$ into non-empty Borel subsets, $(w_e)_{e\in
E}$ is a family of Borel measurable self-maps on the metric space
such that $w_e\left(K_{i(e)}\right)\subset K_{t(e)}$ for all $e\in
E$ and $(p_e)_{e\in E}$ is a family of Borel measurable
probability functions on $K$ (i.e. $p_e(x)\geq 0$ for all $e\in E$
and $\sum_{e\in E}p_e(x)=1$ for all $x\in K$) (associated with the
maps) such that each $p_e$ is zero on the complement of
$K_{i(e)}$.

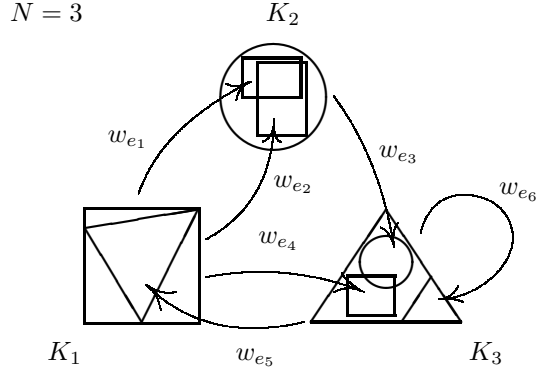
\begin{center}
\unitlength 1mm
\begin{picture}(70,70)\thicklines
\put(35,50){\circle{20}} \put(10,20){\framebox(15,15)}
\put(40,20){\line(2,3){10}} \put(40,20){\line(4,0){20}}
\put(50,35){\line(2,-3){10}} \put(5,15){$K_1$} \put(34,60){$K_2$}
\put(61,15){$K_3$} \put(31,50){\framebox(7.5,5)}
\put(33,45){\framebox(6.25,9.37)} \put(50,28){\circle{7.5}}
\put(45,21){\framebox(6,5)} \put(10,32.5){\line(6,1){15}}
\put(10,32.5){\line(3,-5){7.5}} \put(17.5,20){\line(1,2){7.5}}
\put(52,20){\line(2,3){4}} \put(13,44){$w_{e_1}$}
\put(35,38){$w_{e_2}$} \put(49,42){$w_{e_3}$}
\put(33,30.5){$w_{e_4}$} \put(30,15){$w_{e_5}$}
\put(65,37){$w_{e_6}$} \put(30,5){ Fig. 1} \put(0,60){$N=3$}
\thinlines \linethickness{0.1mm} \bezier{300}(17,37)(20,46)(32,52)
\bezier{50}(32,52)(30.5,51.7)(30,49.5)
\bezier{50}(32,52)(30,51)(28.7,51.7)
\bezier{300}(26,31)(35,36)(35,47)
\bezier{50}(35,47)(35,44.5)(33.5,44)
\bezier{50}(35,47)(35,44)(36,44) \bezier{300}(43,50)(49,42)(51,30)
\bezier{50}(51,30)(50.5,32)(49.2,32.6)
\bezier{50}(51,30)(50.6,32)(51.5,33.2)
\bezier{300}(39,20)(26,17)(18,25)
\bezier{50}(18,25)(19.5,24)(20,21.55)
\bezier{50}(18,25)(20,23.5)(22,24)
\bezier{300}(26,26)(37,28)(47,24) \bezier{50}(47,24)(45,25)(43,24)
\bezier{50}(47,24)(45,25)(44,26.5)
\bezier{100}(54.5,31.9)(56,37.3)(61,36.9)
\bezier{100}(61,36.9)(64.5,36.5)(66,34)
\bezier{100}(66,34)(68,30.5)(64.9,26.8)
\bezier{100}(64.9,26.8)(61.6,23.3)(57,23)
\bezier{50}(57,23)(58.5,23.3)(60.1,22.7)
\bezier{50}(57,23)(58.8,23.3)(59.5,24.8)
\end{picture}
\end{center}

We call a Markov system $\left(K_{i(e)},w_e,p_e\right)_{e\in E}$
{\it contractive} with an {\it average contracting rate} $0<a<1$
iff it satisfies the following {\it condition of contractiveness
on average}
\begin{equation*}
 \sum\limits_{e\in E}p_e(x)d(w_ex,w_ey)\leq ad(x,y)\mbox{ for all
}x,y\in K_i,\ i=1,...,N.
\end{equation*}

Markov system $\left(K_{i(e)},w_e,p_e\right)_{e\in E}$ determines
a Markov operator $U$ on the set of all bounded Borel measurable
functions $\mathcal{L}^0(K)$ by
\[Uf:=\sum\limits_{e\in E}p_ef\circ w_e\mbox{ for all
 }f\in\mathcal{L}^0(K)\] and its adjoint operator $U^*$ on the set of all Borel probability
 measures $P(K)$ by
\[U^*\nu(f):=\int U(f)d\nu\mbox{ for all }f\in\mathcal{L}^0(K)\mbox{ and }\nu\in P(K).\]

 We say a probability measure $\mu$ is an {\it invariant probability measure} of
 the Markov system iff it is a stationary initial distribution of
 the associated Markov process, i.e.
 \[U^*\mu=\mu.\]
A Borel probability measure $\mu$ is called {\it
 attractive} measure of the CMS if
 \[{U^*}^n\nu\stackrel{w^*}{\to}\mu\mbox{ for all }\nu\in P(K),\]
 where $w^*$ means weak$^*$ convergence.
Note that an attractive probability measure is a unique invariant
probability measure of the CMS if $U$ maps continuous functions on
continuous functions. We will denote the space of all bounded
continuous functions by $C_B(K)$.

The main result in \cite{Wer1} concerning the uniqueness of the
invariant measure is the following (see Lemma 1 and Theorem 2 in
\cite{Wer1}).
\begin{theo}\label{Th}
 Suppose $\left(K_{i(e)},w_e,p_e\right)_{e\in E}$ is an irreducible CMS such that  $K_1$,$K_2$,
 ...,$K_N$ partition $K$ into non-empty open subsets, each $p_e|_{K_{i(e)}}$ is
 Dini-continuous  and there exists $\delta>0$ such that $p_e|_{K_{i(e)}}\geq\delta$
 for all $e\in E$. Then:\\
(i) The CMS has a unique invariant Borel probability measure
$\mu$, and $\mu(K_i)>0$ for all $i=1,...,N$.\\
(ii) If in addition the CMS is aperiodic, then
\[U^nf(x)\to\mu(f)\mbox{ for all }x\in K\mbox{ and }f\in C_B(K)
\] and the convergence is uniform on bounded subsets, i.e. $\mu$ is an
attractive probability measure of the CMS.
\end{theo}
A function $h:(X,d)\longrightarrow\mathbb{R}$ is called {\it
Dini-continuous} iff
 for some $c >0$  \[\int_0^c\frac{\phi(t)}{t}dt<\infty\]
  where $\phi$ is {\it the modulus of uniform continuity} of $f$, i.e.
     \[\phi(t):=\sup\{|h(x)-h(y)|:d(x,y)\leq t,\ x,y\in X\}.\]
It is easily seen that the Dini-continuity is weaker than the
H\"{o}lder and stronger than the uniform continuity.  There  is a
well known characterization of the Dini-continuity, which will be
useful later.
\begin{lemma}\label{Dc}
 Let $0<c<1$ and $b>0$.
  A  function $h$ is Dini-continuous  iff
 \[\sum_{n=0}^\infty\phi\left(bc^n\right)<\infty\] where $\phi$ is
 the modulus of uniform continuity of $h$.
\end{lemma}
The proof is simple (e.g. see \cite{Wer1}).

Further, with the Markov system also is associated a measure
preserving transformation $S:(\Sigma,\mathcal{B}(\Sigma),
M)\longrightarrow(\Sigma,\mathcal{B}(\Sigma), M)$, which we call a
{\it generalized Markov shift}, where
$\Sigma:=\{(...,\sigma_{-1},\sigma_0,\sigma_1,...):\sigma_i\in E\
\forall i\in\mathbb{Z}\}$ is the {\it code space} provided with
the product topology, $\mathcal{B}(\Sigma)$ denotes Borel
$\sigma$-algebra on $\Sigma$ and $M$ is a {\it generalized Markov
measure} on $\mathcal{B}(\Sigma)$ given by
\[M\left(_m[e_1,...,e_k]\right):=\int
p_{e_1}(x)p_{e_2}(w_{e_1}x)...p_{e_k}(w_{e_{k-1}}\circ...\circ
w_{e_1}x)d\mu(x)\] for every thin cylinder set
$_m[e_1,...,e_k]:=\{\sigma\in\Sigma:\
\sigma_m=e_1,...,\sigma_{m+k-1}=e_k\}$, $m\in\mathbb{Z}$, where
$\mu$ is an invariant Borel probability measure of the Markov
system, and $S$ is the usual left shift map on $\Sigma$. It is
easy to verify that $S$ preserves measure $M$, since $U^*\mu=\mu$
(see \cite{Wer3}).

For a CMS, this two pictures, Markovian and dynamical, are related
by a {\it coding map} $F:(\Sigma,\mathcal{B}(\Sigma),
M)\longrightarrow K$ which was constructed in \cite{Wer3}. It is
defined, if $K$ is a complete metric space and $p_e|_{K_i(e)}$'s
are Dini-continuous, by
\[F(\sigma):=\lim\limits_{m\to-\infty}w_{\sigma_0}\circ
w_{\sigma_{-1}}\circ...\circ w_{\sigma_m}x_{i(\sigma_m)}\mbox{ for
}M\mbox{-a.e. }\sigma\in\Sigma\] where $x_i\in K_i$ for each
$i=1,...,N$ (the coding map does not depend on the choice of
$x_i$'s modulo an $M$-zero set). This coding map is the key for
our calculation.

\begin{Example}
  Let $G:=(V,E,i,t)$ be a finite irreducible directed
  (multi)graph. Let $\Sigma^-_G:=\{(...,\sigma_{-1},\sigma_0):\
  \sigma_m\in E\mbox{ and }
t(\sigma_m)=i(\sigma_{m-1})\ \forall
m\in\mathbb{Z}\setminus\mathbb{N}\}$ ({\it one-sided
  subshift of finite type} associated with $G$) endowed with the
  metric $d(\sigma,\sigma'):=2^k$ where $k$ is the smallest
  integer with $\sigma_i=\sigma'_i$ for all $k<i\leq 0$. Let $g$
  be a positive, Dini-continuous function on $\Sigma_G$ such that
  \[\sum\limits_{y\in T^{-1}(\{x\})}g(y)=1\mbox{ for all }x\in\Sigma_G\]
  where $T$ is the right shift map on $\Sigma^-_G$. Set
  $K_i:=\left\{\sigma\in\Sigma^-_G:t(\sigma_0)=i\right\}$  for
  every $i\in V$ and, for
  every $e\in E$,
  \[w_e(\sigma):=(...,\sigma_{-1},\sigma_{0},e),\ p_e(\sigma):=g(...,\sigma_{-1},\sigma_{0},e)
  \mbox{ for all }\sigma\in K_{i(e)}.\]
  Obviously, maps $(w_e)_{e\in E}$ are contractions. Therefore,
  $\left(K_{i(e)}, w_e, p_e\right)_{e\in E}$
   defines a CMS. An invariant probability measure of this CMS is
   called a $g$-measure. This notion was introduced by M. Keane
   \cite {Ke}. The Kolmogorov-Sinai entropy of the generalized Markov shift, which is a natural
   extension of the $g$-measure $\mu$ in this case, is well known:
   \[h_M(S)=\sum\limits_{e\in E}\int g(e\sigma)\log g(e\sigma)d\mu(\sigma),\]
   where $e\sigma:=(...,\sigma_{-1},\sigma_0,e)$
   (see \cite{Le}, \cite{W1}).
\end{Example}

So, we are going to extend this result to a more general CMS, e.g.
as in the next example.

\begin{Example}
 Let $\mathbb{R}^2$ be normed by $\|.\|_1$. Let
 $K_1:=\{(x,y)\in\mathbb{R}^2:\ y\geq 1/2\}$ and $K_2:=\{(x,y)\in\mathbb{R}^2:\ y\leq
 -1/2\}$. Consider the following maps on $\mathbb{R}^2$:
 \begin{eqnarray*}
   &&w_1\left(\begin{array}{c}x\\ y\end{array}\right):=\left(\begin{array}{c}
   -\frac{1}{2}x-1\\-\frac{3}{2}y+\frac{1}{4}\end{array}\right),\
   w_2\left(\begin{array}{c}x\\ y\end{array}\right):=\left(\begin{array}{c}
   -\frac{3}{2}x+1\\\frac{1}{4}y+\frac{3}{8}\end{array}\right),\\
  &&w_3\left(\begin{array}{c}x\\ y\end{array}\right):=\left(\begin{array}{c}
   -\frac{1}{2}|x|+1\\-\frac{3}{2}y-\frac{1}{4}\end{array}\right),\mbox{
   and }
   w_4\left(\begin{array}{c}x\\ y\end{array}\right):=\left(\begin{array}{c}
   \frac{3}{2}|x|-1\\-\frac{1}{4}y+\frac{3}{8}\end{array}\right)
 \end{eqnarray*}
 with probability functions
 \begin{eqnarray*}
  && p_1\left(\begin{array}{c}x\\ y\end{array}\right):=\left(\frac{1}{15}\sin^2\|(x,y)\|_1+\frac{53}{105}
   \right)1_{K_1}(x,y),\\
  && p_2\left(\begin{array}{c}x\\ y\end{array}\right)
   :=\left(\frac{1}{15}\cos^2\|(x,y)\|_1+\frac{3}{7}\right)
   1_{K_1}(x,y),\\
   && p_3\left(\begin{array}{c}x\\ y\end{array}\right):=\left(\frac{1}{15}\sin^2\|(x,y)\|_1+\frac{53}{105}
   \right)1_{K_2}(x,y)\mbox{ and }\\
  && p_4\left(\begin{array}{c}x\\ y\end{array}\right)
   :=\left(\frac{1}{15}\cos^2\|(x,y)\|_1+\frac{3}{7}\right) 1_{K_2}(x,y)  .
 \end{eqnarray*}
 A simple calculation shows that $(K_{i(e)},w_e,p_e)_{e\in\{1,2,3,4\}}$, where $i(1)=i(2)=1$ and $i(3)=i(4)=2$,
 defines a CMS with an average
 contracting rate $209/210$. Note that none of the maps are contractive (By Theorem 2 in \cite{Wer1}, it
 has a unique (attractive) invariant probability measure.)
\end{Example}

As this paper had been known for a year as a preprint, I was
informed by Wojciech Slomczynski that a similar entropy formula for
some partial iterated function systems plays a central role in his
recent book \cite{S}. However, the result presented in this paper is
not covered by any results from his book. Moreover, our computation
is based on the existence of the coding map for some contractive
Markov systems \cite{Wer3}. The latter is a nontrivial fact which
seems to be of a fundamental nature in relation to the uniqueness of
the invariant probability measure for a contractive Markov system.

\section{Main Part}
Let $\left(K_{i(e)},w_e,p_e\right)_{e\in E}$ be a contractive
Markov system with the average contracting rate $0<a<1$ and an
invariant Borel probability measure $\mu$. We assume that: $(K,d)$
is a metric space in which sets of finite diameter are relatively
compact and the family $K_1,...,K_N$ partitions $K$ into non-empty
open subsets; each probability function $p_e|_{K_{i(e)}}$ is
Dini-continuous and bounded away from zero by $\delta>0$; the set
of edges $E$ is finite and the map $i:E\longrightarrow V$ is
surjective. Note that the assumption on the metric space implies
that it is locally compact separable and complete.

First, we prove what seems to be the main lemma for the generalized
Markov shift  associated with a contractive Markov system. This
lemma is also used in \cite{Wer5}. For that we need to define some
measures on the product space $K\times\Sigma$.

Denote by $\mathcal{A}$  the finite $\sigma$-algebra generated by
the partition $\{_0[e]:e\in E\}$ of $\Sigma$ and define, for each
integer $m\leq 1$,
\[\mathcal{A}_m:=\bigvee\limits_{i=m}^{+\infty} S^{-i}\mathcal{A},\]
which is the smallest $\sigma$-algebra containing  all finite
 $\sigma$-algebras $\bigvee_{i=m}^{n}
S^{-i}\mathcal{A}$, $n\geq m$. Let $x\in K$. For every integer
$m\leq 1$, let $P_x^m$ be a probability measure on
$\sigma$-algebra $\mathcal{A}_m$ given by
\[P^m_x( _{m}[e_{m},...,e_n])=p_{e_{m}}(x)p_{e_{m+1}}(w_{e_{m}}(x))...p_{e_n}(w_{e_{n-1}}\circ...\circ
w_{e_{m}}(x))\] for all thin cylinders $_{m}[e_{m},...,e_n]$,
$n\geq{m}$. By Lemma 1 from \cite{Wer3}, $x\longmapsto P_x^m(Q)$ is
a Borel measurable function on $K$. Therefore, we can define, for
every integer $m\leq 0$,
\[\tilde M_m\left(A\times Q\right):=\int\limits_{A}P^m_x\left(Q\right)d\mu(x)\]
for all $A\in\mathcal{B}(K)$ and all $Q\in\mathcal{A}_m$. Then
$\tilde M_m$ extends uniquely to a probability measure on the
product $\sigma$-algebra $\mathcal{B}( K)\otimes\mathcal{A}_m$
with
\[\tilde M_m(\Omega)=\int
P^m_x\left(\left\{\sigma\in\Sigma:(x,\sigma)\in\Omega\right\}\right)d\mu(x)\]
for all $\Omega\in\mathcal{B}( K)\otimes\mathcal{A}_m$. Note that
the set of all $\Omega\in\mathcal{B}( K)\otimes\mathcal{A}_m$ for
which the integrand in the above is measurable forms a Dynkin
system which contains the set all rectangles   $A\times Q$,
$A\in\mathcal{B}(K)$, $Q\in\mathcal{A}_m$. As the latter is
$\cap$-stable and generates $\mathcal{B}(K)\otimes\mathcal{A}_m$,
the integrand is measurable for all
$\Omega\in\mathcal{B}(K)\otimes\mathcal{A}_m$. Further, note that
$P^m_x\left(\left\{\sigma\in\Sigma:(x,\sigma)\in\Omega\right\}\right)=\int
1_{\Omega}(x,\sigma)dP_x^m(\sigma)$ for all
$\Omega\in\mathcal{B}(K)\otimes\mathcal{A}_m$. Therefore
\[\int s d\tilde M_m=\int\int s(x,\sigma)dP^m_x(\sigma)d\mu(x)\] for
all $\mathcal{B}(K)\otimes\mathcal{A}_m$-simple functions $s$.
Now, let $\psi$ be a
$\mathcal{B}(K)\otimes\mathcal{A}_m$-measurable and $\tilde
M_m$-integrable function on $K\times\Sigma$. Then the
 usual monotone approximation of positive and negative parts of $\psi$ by simple functions
 and the B. Levi Theorem imply that
\[\int \psi d\tilde M_m=\int\int \psi(x,\sigma)dP^m_x(\sigma)d\mu(x).\]

\begin{lemma}\label{mlMs} Suppose
$C:=\sum_{i=1}^N\int_{K_i}d(x,x_i)d\mu(x)<\infty$ for some $x_i\in
K_i,\ i=1,...,N$. Let $_1[e_1,...,e_n]\subset\Sigma$ be a thin
cylinder set . Let
  $\mathcal{F}:=\bigvee_{i=0}^\infty S^i\mathcal{A}$. Then
 \[E_M\left(1_{_1[e_1,...,e_n]}|\mathcal{F}\right)(\sigma)= P^1_{F(\sigma)}\left(_1[e_1,...,e_n]\right)
 \mbox{ for }M\mbox{-a.e. }\sigma\in\Sigma,\] where $E_M(.|.)$
 denotes conditional expectation with respect to measure $M$.
\end{lemma}
{\it Proof.} We can obviously assume that $(e_1,...,e_n)$ is a path
of the directed graph. Set $\mathcal{F}_m:=\bigvee_{i=0}^m
S^i\mathcal{A}$ for all $m\in\mathbb{Z}\setminus\mathbb{N}$. We
denote  further a $(-m+1)$-tuple by $(\sigma_m,...,\sigma_0)^*$ if
$M\left(_m[\sigma_m,...,\sigma_0]\right)>0 $ (i.e.
$(\sigma_m,...,\sigma_0)$ is a path). Then obviously
\[E_M\left(1_{_1[e_1,...,e_n]}|\mathcal{F}_m\right)(\tilde\sigma)=\sum\limits_{(\sigma_m,...,\sigma_0)^*}
  \frac{\int\limits_{_m[\sigma_m,...,\sigma_0]}
1_{_1[e_1,...,e_n]}dM}{M\left(_m[\sigma_m,...,\sigma_0]\right)}1_{_m[\sigma_m,...,\sigma_0]}(\tilde\sigma)\]
for $M$-a.e. $\tilde\sigma\in\Sigma$. Since
$\left(\mathcal{F}_m\right)_{m\leq 0}$ is an increasing sequence of
$\sigma$-algebras and $\mathcal{F}$ is the smallest $\sigma$-algebra
containing all $\mathcal{F}_m$, it follows by Doob's Martingale
Theorem (e.g. see \cite{Doob} p. 199) that
\begin{eqnarray}\label{mc}
E_M\left(1_{_1[e_1,...,e_n]}|\mathcal{F}_m\right)\to
E_M\left(1_{_1[e_1,...,e_n]}|\mathcal{F}\right)
\end{eqnarray}   $M$-a.e..

Now,  set
\[Z^x_m(\sigma):=w_{\sigma_0}\circ...\circ w_{\sigma_m}x\mbox{
and }Y_m(\sigma):=w_{\sigma_0}\circ...\circ
w_{\sigma_m}x_{i(\sigma_m)}\]  for all $x\in K$, $\sigma\in\Sigma$
and $m\leq 0$. Further, define
\[X_m(\tilde\sigma):=\sum\limits_{(\sigma_m,...,\sigma_0)^*}
  \frac{\int\limits_{K\times _m[\sigma_m,...,\sigma_0]}
d\left(Z^x_m(\bar\sigma),Y_m(\bar\sigma)\right)d\tilde
M_m(x,\bar\sigma)}{\tilde M_m\left(K\times
_m[\sigma_m,...,\sigma_0]\right)}1_{_m[\sigma_m,...,\sigma_0]}(\tilde\sigma)\]
for all $\tilde\sigma\in\Sigma$. Then
\begin{eqnarray*}
 &&\int X_mdM=\int d\left(Z^x_m(\sigma),Y_m(\sigma)\right)d\tilde M_m(x,\sigma)\\
&=&\int\sum\limits_{\sigma_{m},...,\sigma_0}p_{\sigma_{m}}(x)...p_{\sigma_{0}}(w_{\sigma_{-1}}...
w_{\sigma_{m}}x) d(w_{\sigma_{0}}...
w_{\sigma_{m}}x,w_{\sigma_{0}}...
w_{\sigma_{m}}x_{i(\sigma_{m})})d\mu(x)\\
& \leq&
a^{-m+1}\sum\limits_{i=1}^N\int\limits_{K_i}d(x,x_i)d\mu(x)=
a^{-m+1}C.
\end{eqnarray*}
 Set
\[\Omega_{m}:=\left\{\sigma\in \Sigma:\
X_m>a^{\frac{-m+1}{2}} C\right\}\] for all $m\leq 0$ and
$\Omega:=\bigcap_{n\leq 0}\bigcup_{m\leq n}\Omega_m$. Then, by
Markov inequality,
 \[ M\left(\Omega_m\right)\leq
a^{\frac{-m+1}{2}}.\] Hence, $M(\Omega)=0$ and
\[X_m(\sigma)\to 0\mbox{ for all }\sigma\in\Sigma\setminus\Omega.\]

 Now, for $\sigma\in\Sigma$ with
$M\left(_m[\sigma_m,...,\sigma_0]\right)>0$,
\begin{eqnarray}\label{eq}
 & &\left| \frac{\int\limits_{_m[\sigma_m,...,\sigma_0]}
1_{_1[e_1,...,e_n]}dM}{M\left(_m[\sigma_m,...,\sigma_0]\right)}-P^1_{F(\sigma)}\left(_1[e_1,...,e_n]\right)\right|
\nonumber\\
&=&\left|\frac{\int p_{\sigma_m}(x)...
p_{\sigma_0}(w_{\sigma_{-1}}... w_{\sigma_m}x)
p_{e_1}(Z^x_m(\sigma))...p_{e_n}(w_{e_{n-1}}... w_{e_1}
Z^x_m(\sigma))d\mu(x)} {\int p_{\sigma_m}(x)...
p_{\sigma_0}(w_{\sigma_{-1}}... w_{\sigma_m}x)d\mu(x)}\right.\nonumber\\
& & \left.- p_{e_1}(F(\sigma))...p_{e_n}(w_{e_{n-1}}... w_{e_1}
F(\sigma))\right|.
\end{eqnarray}
 Set
$p(x):=p_{e_1}(x)...p_{e_n}(w_{e_{n-1}}... w_{e_1} x)$, $x\in K$.
Note that the average contractiveness condition and the
boundedness away from zero of the probability functions on their
vertex sets imply that each map $w_e|_{K_{i(e)}}$ is Lipschitz.
Since each $p_e|_{K_{i(e)}}$ is Dini-continuous, it follows that
each function $p_{e_k}\circ w_{e_{k-1}}\circ...\circ
w_{e_1}|_{K_{i(e_1)}}$, $1\leq k\leq n$, is Dini-continuous. As
bounded Dini-continuous functions form an algebra,
$p|_{K_{i(e_1)}}$ is also Dini-continuous. Let $\eta$ be the
modulus of uniform continuity of $p|_{K_{i(e_1)}}$. By the
Sublemma from \cite{BDEG2}, there exists
$\psi:[0,\infty)\longrightarrow[0,\infty)$ such that
$\psi(t)\geq\eta(t)$ for all $t$, $\psi(t)/t$ is non-increasing,
and $\int_0^1\psi(t)/tdt<\infty$. Set
\[\beta(u):=\frac{1}{1-a}\int\limits_0^{ua^{-1}}\frac{\psi(t)}{t}dt\mbox{ for }u\geq
0.\] Then $\beta$ is continuous, concave and $\beta(0)=0$.
Moreover,
\[\beta(u)\geq\frac{1}{1-a}\int\limits_u^{ua^{-1}}\frac{\psi(t)}{t}dt\geq
\frac{1}{1-a}\frac{\psi(u)}{ua^{-1}}u(a^{-1}-1)=\psi(u)\geq
\eta(u)\] for all $u> 0$. Hence, $\eta(u)\leq\beta(u)$ for all
$u\geq 0$. Therefore,
\begin{eqnarray*}
 (\ref{eq})&\leq&\left|\frac{\int\limits_{K\times\
_m[\sigma_m,...,\sigma_0]} p\circ Z^x_m(\bar\sigma)d\tilde
M_m(x,\bar\sigma)} {\tilde M_m\left(K\times\
_m[\sigma_m,...,\sigma_0]\right)}-p\circ Y_m(\sigma)\right|+
\left|p\circ Y_m(\sigma)-p\circ F(\sigma)\right|\\
 &\leq&\frac{\int\limits_{K\times\
_m[\sigma_m,...,\sigma_0]} \left|p\circ Z^x_m(\bar\sigma)-p\circ
Y_m(\bar\sigma)\right|d\tilde M_m(x,\bar\sigma)} {\tilde
M_m\left(K\times\ _m[\sigma_m,...,\sigma_0]\right)}+
\left|p\circ Y_m(\sigma)-p\circ F(\sigma)\right|\\
&\leq&\frac{\int\limits_{K\times\ _m[\sigma_m,...,\sigma_0]}
\beta\left(d\left(Z^x_m(\bar\sigma),Y_m(\bar\sigma)\right)\right)d\tilde
M_m(x,\bar\sigma)} {\tilde M_m\left(K\times\
_m[\sigma_m,...,\sigma_0]\right)}+
\left|p\circ Y_m(\sigma)-p\circ F(\sigma)\right|\\
&\leq&\beta\left(\frac{\int\limits_{K\times\
_m[\sigma_m,...,\sigma_0]}
d\left(Z^x_m(\bar\sigma),Y_m(\bar\sigma)\right)d\tilde
M_m(x,\bar\sigma)} {\tilde M_m\left(K\times\
_m[\sigma_m,...,\sigma_0]\right)}\right)+ \left|p\circ
Y_m(\sigma)-p\circ F(\sigma)\right|\\
&=& \beta\circ X_m(\sigma)+\left|p\circ Y_m(\sigma)-p\circ
F(\sigma)\right|.
\end{eqnarray*}
Hence
\begin{eqnarray*}
\left|E_M\left(1_{_1[e_1,...,e_n]}|\mathcal{F}_m\right)(\sigma)-P^1_{F(\sigma)}\left(_1[e_1,...,e_n]\right)\right|
\leq\beta\circ X_m(\sigma)+\left|p\circ Y_m(\sigma)-p\circ
F(\sigma)\right|
\end{eqnarray*}
for $M$-a.e. $\sigma\in\Sigma$. By Corollary 1 from \cite{Wer3}
and the continuity of $p$ on $K$, the second term also converges
to zero $M$-a.e.. Thus
\[\left|E_M\left(1_{_1[e_1,...,e_n]}|\mathcal{F}_m\right)(\sigma)-P^1_{F(\sigma)}\left(_1[e_1,...,e_n]\right)\right|
\to 0\mbox{ as }m\to-\infty\] for $M$-a.e. $\sigma\in\Sigma$. With
(\ref{mc}), this implies that
\[E_M\left(1_{_1[e_1,...,e_n]}|\mathcal{F}\right)(\sigma)=P^1_{F(\sigma)}\left(_1[e_1,...,e_n]\right)\]
for $M$-a.e. $\sigma\in\Sigma$. \hfill$\Box$

\begin{prop}\label{iMm}
If invariant probability measure $\mu$ is unique, then
\[F(M)=\mu.\]
\end{prop}
{\it Proof.}\noindent Let $U^*$ be the adjoint of the Markov
operator associated with the CMS. It is sufficient to show that
$U^*F(M)=F(M)$, since $\mu$ is the unique invariant probability
measure. Let $f\in C_B(K)$. Then
\begin{eqnarray*}
 U^*F(M)(f)&=&\sum\limits_{e\in E}\int p_ef\circ w_edF(M)=\sum\limits_{e\in E}\int p_e\circ Ff\circ w_e
 \circ FdM.
\end{eqnarray*}
 Let $e\in E$.  By  Theorem 1 (iv) from \cite{Wer1},
 $\sum_{i=1}^N\int_{K_i}d(x,x_i)d\mu(x)<\infty$.
 Therefore, by Lemma \ref{mlMs},
\[E_M\left(\left.1_{_1[e]}\right|\mathcal{F}\right)(\sigma)=P^1_{F(\sigma)}\left(
_1[e]\right)=p_e\circ F(\sigma)\] for $M$-a.e. $\sigma\in\Sigma$.
Since $f\circ w_e \circ F$ is bounded and
$\mathcal{F}$-measurable, it follows by a well known property of
the conditional expectation that
\[E_M\left(\left.1_{_1[e]}f\circ w_e \circ F\right|\mathcal{F}\right)(\sigma)=
p_e\circ F(\sigma)f\circ w_e \circ F(\sigma)\] for $M$-a.e.
$\sigma\in\Sigma$. Hence, by the shift invariance of $M$,
\begin{eqnarray*}
  U^*F(M)(f)&=&\sum\limits_{e\in E}\int  1_{_1[e]}(\sigma)f\circ w_e
 \circ F(\sigma)dM(\sigma)\\
 &=&\sum\limits_{e\in E}\int  1_{_1[e]}\circ S^{-1}(\sigma)f\circ w_e
 \circ F\circ S^{-1}(\sigma)dM(\sigma)\\
 &=&\sum\limits_{e\in E}\int  1_{_0[e]}(\sigma)f
 \circ F(\sigma)dM(\sigma)\\
 &=&F(M)(f).
\end{eqnarray*}
\hfill$\Box$

\begin{theo}
 Let $h_M(S)$ be Kolmogorov-Sinai entropy of the
generalized Markov shift associated with the contractive Markov
system. \\
$(i)$ If $\sum_{i=1}^N\int_{K_i}d(x,x_i)d\mu(x)<\infty$ for some
$x_i\in K_i$, $i=1,...,N$, then
  \[h_M(S)=-\sum\limits_{e\in E}\int\limits_{K_{i(e)}} p_e\log
  p_edF(M).\]
$(ii)$ If invariant probability measure $\mu$ is unique, then
\[h_M(S)=-\sum\limits_{e\in E}\int\limits_{K_{i(e)}} p_e\log
  p_ed\mu.\]
\end{theo}
{\it Proof.} It is well known that $h_M(S)=h_M(S^{-1})$ (e.g.
Theorem 4.13 in \cite {W}) and, by the Kolmogorov-Sinai Theorem
(e.g. Theorem 4.17 in \cite{W}),
$h_M(S^{-1})=h_M(S^{-1},\mathcal{A})$. Further, using the notion of
conditional entropy (e.g. Theorem 4.3 (ix) and Theorem 4.14 in
\cite{W}),
\[h_M(S^{-1},\mathcal{A})=\mathcal H\left(S^{-1}\mathcal A\left/\bigvee\limits_{i=1}^\infty S^{i-1}\mathcal
A\right.\right).\] Set $\mathcal F:=\bigvee_{i=0}^\infty
S^i\mathcal A$. Hence
\begin{eqnarray*}
 h_M(S)=-\sum\limits_{e\in E}\int E\left(\left.1_{_1[e]}\right|\mathcal F\right)\log E\left(\left.1_{_1[e]}\right|
\mathcal F\right)dM.
\end{eqnarray*}
 By the assumption in $(i)$, Lemma \ref{mlMs} implies that
\begin{eqnarray*}
  E\left(\left.1_{_1[e]}\right|\mathcal F\right) = p_e\circ
  F\ M\mbox{-a.e.}
\end{eqnarray*} for each $e\in E$.
Hence, with $0\log 0=0$,  we have
\begin{eqnarray*}
h_M(S)=-\sum\limits_{e\in E}\int p_e\circ F\log\left( p_e\circ
F\right) dM=-\sum\limits_{e\in E}\int\limits_{K_{i(e)}} p_e\log
p_e dF(M).
\end{eqnarray*}

By the assumption in $(ii)$, Theorem 1(iv) in \cite{Wer1} and
Proposition \ref{iMm} imply that
\[h_M(S)=-\sum\limits_{e\in E}\int\limits_{K_{i(e)}} p_e\log p_e d\mu.\]
\hfill$\Box$

\subsection*{Acknowledgements} This work was done as a part of my PhD thesis at the University of St Andrews.
I would like to thank: EPSRC and School of Mathematics
 and Statistics of University of St Andrews for providing me with a scholarship and excellent working
 conditions in St Andrews, Professor K. J. Falconer for valuable comments on the first draft of the paper,
 my supervisor Lars Olsen for his interest
 in my work, his support and many fruitful discussions.

\end{document}